\documentclass[10pt,dvipsnames,reqno]{amsart}

\pdfoutput=1

\usepackage{filecontents}

\usepackage[style=numeric, backend=bibtex, backref=true, bibencoding=utf8, datamodel=mrnumber, maxbibnames=6, sorting=nyt]{biblatex}

\DeclareFieldFormat{mrnumber}{%
  \ifhyperref
    {\href{http://www.ams.org/mathscinet-getitem?mr=1#1}{MR#1}}
    {MR#1}}

\DeclareFieldFormat{mrclass}{#1}

\DeclareNameAlias{bymrreviewer}{byeditor}

\newbibmacro*{mrinfo}{%
  \printfield{mrnumber}%
  \iffieldundef{mrclass}
    {\setunit*{\addcomma\space}}
    {\setunit*{\addspace}}%
  \printfield{mrclass}%
  \setunit*{\addcomma\space}%
  \ifnameundef{mrreviewer}
    {}
    {\bibstring{byreviewer}%
     \setunit{\addspace}%
     \printnames[bymrreviewer]{mrreviewer}}}

\newtoggle{bbx:mrinfo}
\DeclareBibliographyOption[boolean]{mrinfo}[true]{\settoggle{bbx:mrinfo}{#1}}
\renewbibmacro*{doi+eprint+url}{%
  \iftoggle{bbx:doi}
    {\printfield{doi}}
    {}%
  \newunit\newblock
  \iftoggle{bbx:mrinfo}
    {\usebibmacro{mrinfo}}
    {}%
  \newunit\newblock
  \iftoggle{bbx:eprint}
    {\usebibmacro{eprint}}
    {}%
    }

\usepackage{adjustbox, afterpage, amsbsy, amsmath, amssymb, amsthm,
  caption, diagbox, enumitem, fixmath, graphicx, mathtools,
  mdframed, pigpen, pdfpages, pgfkeys, soul, subcaption,
  xspace, url}

\usepackage{symbols}
\usepackage[utf8]{inputenc}
\usepackage[nonewpage]{imakeidx}
\usepackage{interval} %%% formats intervals in posets
\intervalconfig{soft open fences}
\usepackage{xcolor}
\usepackage{hyperref}
\hypersetup{
     colorlinks=true,
     linkcolor=blue,
     filecolor=blue,
     citecolor=blue,
     urlcolor=blue,
   }

\allowdisplaybreaks
\usepackage[text={13.8cm,25.2cm}, centering]{geometry}

\makeatletter
\def\msection{\@startsection{section} %name
  {1} % level
  {0pt} % indent
  {-3.5ex plus -1ex minus -0.2ex} % beforeskip
  {-2.3ex plus -.2ex} % afterskip
  {\normalfont\large\bfseries} % style
}

\makeatother

\newcommand{\mc}[1]{\mbox{$\mathfrak{C}_{{#1}}$}\xspace}

\theoremstyle{plain}

\theoremstyle{definition}

\newtheorem{rem}{Remark}[section]

\title[Finite Sum of Subobjects]{Internal Neighbourhood Structures III: Finite Sum of Subobjects}
\author{Partha Pratim Ghosh}
\address{Department of Mathematical Sciences \\
  University of South Africa \\
  Unisa Science Campus \\
  corner of Christiaan de Wet \& Pioneer Avenue \\
  Florida 1709 \\
  Johannesburg, Gauteng \\
  South Africa}
\email{ghoshpp@unisa.ac.za}

\begin{document}

%%% need to know the dimensions for this document

% \convertto{cm}{1em} \convertto{cm}{1ex}
% \convertto{pt}{1em} \convertto{pt}{1ex}
% \large \convertto{cm}{1em} \convertto{cm}{1ex}
% \large \convertto{pt}{1em} \convertto{pt}{1ex}

% \Large \convertto{cm}{1em} \convertto{cm}{1ex}
% \Large \convertto{pt}{1em} \convertto{pt}{1ex}

% \Huge \convertto{cm}{1em} \convertto{cm}{1ex}
% \Huge \convertto{pt}{1em} \convertto{pt}{1ex}

% \small \convertto{cm}{1em} \convertto{cm}{1ex}
% \small \convertto{pt}{1em} \convertto{pt}{1ex}

% \scriptsize \convertto{cm}{1em} \convertto{cm}{1ex}
% \scriptsize \convertto{pt}{1em} \convertto{pt}{1ex}

%%%%%%%%%%%%%%%%%%%%%%%%%%%%%%%%%%%%%%%%%%%%%%

%%% avoid big vertical spaces between equation environments
\setlength{\abovedisplayskip}{0pt}
\setlength{\abovedisplayshortskip}{0pt}
\setlength{\belowdisplayskip}{0pt}
\setlength{\belowdisplayshortskip}{0pt}
%%%%%%%%%%%%%%%%%%%%%%%%%%%%%%%%%%%%%%%%%%%%%%%%%%

\begin{abstract}
  The notion of an internal preneighbourhood space on a finitely complete category with finite coproducts and a proper \fact{\textsf{E}}{\textsf{M}} system such that for each object $X$ the set of \textsf{M}-subobjects of $X$ is a complete lattice was initiated in \cite[][]{2020}. The notion of a closure operator, closed morphism and its near allies investigated in \cite[][]{2021-clos}{.} The present paper provides structural conditions on the triplet $(\Bb{A}, \mathsf{E}, \mathsf{M})$ (with \Bb{A}{} lextensive) equivalent to the set of \textsf{M}-subobjects of an object closed under finite sums. Equivalent conditions for the set of closed embeddings (closed morphisms) closed under finite sums is also provided. In case when lattices of admissible subobjects (respectively, closed embeddings) are closed under finite sums, the join semilattice of admissible subobjects (respectively, closed embeddings) of a finite sum is shown to be a biproduct of the component join semilattices. Finally, it is shown whenever the set of closed morphisms is closed under finite sums, the set of proper (respectively, separated) morphisms are also closed under finite sums. This leads to equivalent conditions for the full subcategory of compact (respectively, Hausdorff) preneighbourhood spaces to be closed under finite sums.
\end{abstract}



\maketitle

\msection{Introduction}
\label{sec:introduction}

The notion of an internal preneighbourhood space in a context (used to
mean a triplet $\CAL{A} = (\Bb{A}, \mathsf{E}, \mathsf{M})$, where
\Bb{A} is a finitely complete category with finite coproducts and a
proper \fact{\textsf{E}}{\textsf{M}} structure such that for each
object $X$, the set \Sub{X}{\mathsf{M}}{} of \textsf{M}-subobjects of
$X$, also called \emph{admissible subobjects}, is a complete lattice)
appeared in \cite{2020}{.} In \cite[][]{2021-clos} a closure operator
\Arr{\cls{}{\mu}}{\Sub{X}{\mathsf{M}}}{\Sub{X}{\mathsf{M}}}{} on the
lattice \Sub{X}{\textsf{M}} of admissible subobjects is
introduced. Along with near allies of closed morphisms --- closed
morphisms, dense morphisms, proper morphisms, separated morphisms,
perfect morphisms are also investigated in detail. Special proper
(respectively, separated) morphisms give rise to compact
(respectively, Hausdorff) internal preneighbourhood spaces, are also
introduced \cite[see][for details]{2021-clos}.  The purpose of the
present paper is to find conditions under which finite sum of closed
morphisms is closed. This property is necessary to ensure, in particular, the sum of compact or Hausdorff preneighbourhood spaces to be compact or Hausdorff again. Towards this goal, the base category \Bb{A}{} of
the context \CAL{A}{} is restricted to a lextensive category
\cite[see][for details]{CarboniLackWalters1993}, wherein the sums are
well behaved. The main results of this paper are:
\begin{enumerate}[label=(\alph*)]
\item In an extensive context finite sum of admissible subobjects is
  again an admissible subobject if and only if the monomorphisms in
  \textsf{E} between finite sums are stable under pullbacks along
  coproduct injections (Theorem
  \ref{Extensivity>pbstabmonicformalepis<>admsumclosed}{}).
    
\item Theorem \ref{sumofclosedsubobjects}{} shows in an extensive
  context the following are equivalent:
  \begin{enumerate}[label=(\Alph*)]
  \item Finite sum of closed embeddings (i.e., closed morphisms which
    are admissible monomorphisms) is a closed embedding.
    
  \item Finite sum of admissible subobjects is an admissible subobject
    and each coproduct injection is a closed embedding.
    
  \item Dense morphisms between finite sums is stable under pullbacks
    along coproduct injections.
  \end{enumerate}
  
  \noindent{}
  In extensive contexts with admissible subobjects closed under finite
  sums, closed morphisms are closed under finite sums if and only if
  closed embeddings are closed under finite sums (Corollary
  \ref{cor:sum-of-closed-morphisms}{}). Furthermore, if coproduct
  injections are closed, then the closed \textsf{E}-monomorphisms
  between finite sums is stable under pullbacks along coproduct
  injections (Theorem \ref{pbstabilityofclosedmonicformalepi}{}).
  
\item In extensive context, when finite sum of admissible subobjects
  (respectively, closed embeddings) is an admissible subobject
  (respectively, closed embedding) then the join semilattice of admissible
  subobjects (respectively, closed embeddings) of finite sum is a
  biproduct of the corresponding component join semilattices (Theorem
  \ref{subobjectofsumisbiproduct}{}).

\item Theorem \ref{thm:sum-proper} shows in an extensive context,
  finite sum of proper morphisms is a proper morphism
  \cite[see][Definition 6]{2021-clos}. In particular, finite sum of
  compact preneighbourhood spaces \cite[see][Definition
  6.2]{2021-clos} is compact if and only if for each internal
  preneighbourhood space \opair{X}{\mu}{} the unique morphism
  \Arr{\inita{X}}{\inito}{X} and the \emph{codiagonal} morphism
  \Arr{\Opair{\id{X}}{\id{X}}}{X+X}{X}{} are both proper morphism ---
  the colimiting objects \inito{} and $X + X$ being given the largest
  preneighbourhood neighbourhood systems granted by the topologicity
  of the forgetful functor \Arr{U}{\pNHD{\Bb{A}}}{\Bb{A}} from the
  category of internal preneighbourhood spaces and preneighbourhood
  morphisms on \Bb{A}{} to \Bb{A}{} \cite[see][Theorem 4.8(a)]{2020}.
 
\item Theorem \ref{thm:sum-of-separated} shows in an extensive
  context, finite sum of separated morphisms is a separated morphism
  \cite[see][Definition 7]{2021-clos}. In particular, finite sum of
  Hausdorff preneighbourhood spaces \cite[see][Definition
  7.2]{2021-clos} is Hausdorff if and only if $\termo+\termo$ (with
  the largest preneighbourhood system granted by the topologicity of
  \Arr{U}{\pNHD{\Bb{A}}}{\Bb{A}}{}) is Hausdorff.

\end{enumerate}

The notation and terminology adopted in this paper are largely in line
with the usage in \cite{Maclane1998} or \cite{Borceux1994_01}.  Apart
from this, some specific notations and terms are explained here. Given
the proper \fact{\textsf{E}}{\textsf{M}} system, the morphisms of
\textsf{E} are depicted with arrows like \fepi{}{}{}{} while the
morphisms of \textsf{M} are depicted with arrows like
$\xymatrix{ {} \ar@{>->}[r] & {}}$.  If \Arr{f}{X}{Y}{} be a morphism,
$m \in \Sub{X}{\mathsf{M}}$, $n \in \Sub{Y}{\mathsf{M}}$ then the image of
$m$ (respectively, preimage of $n$) under $f$ is \img{f}{m}
(respectively, \finv{f}{n}{}){,} where
$\comp{f}{m} = \comp{(\img{f}{m})}{\rest{f}{M}}$ (respectively,
$\comp{f}{(\finv{f}{n})} = \comp{n}{f_n}$) is the
\fact{\textsf{E}}{\textsf{M}} of \comp{f}{m}{} (respectively, pullback
of $n$ along $f$), \rest{f}{M}{} is the \emph{restriction of $f$ on
  $m$} (respectively, $f_{n}$ is the \emph{corestriction} of $f$ on
$n$).  Furthermore, for each object $X$, the unique morphism
\Arr{\inita{X}}{\inito{}}{X} from the initial object \inito{} has the
\fact{\textsf{E}}{\textsf{M}}
$\inita{X} = \comp{\sigma_X}{\inita{\emptyset_X}}$ as depicted by the diagram
$\xymatrix{ {\inito} \ar@{->>}[r]^-{\inita{\emptyset_X}} & {\emptyset_X}
  \ar@{>->}[]!<12pt,0pt>;[r]^-{\sigma_{X}} & {X} }$ making
$\emptyset_X \in \Sub{X}{\mathsf{M}}$ the smallest subobject of $X$. Finally,
\PosZero{} is the category with objects partially ordered sets having
smallest element, morphisms order preserving maps preserving the
smallest element, and in a pointed category \Bb{X}, the zero object is
denoted by \zeroo{},
$\xymatrix{ {A} \ar[r]^-{\zeroaf{A}} & {\zeroo} \ar[r]^-{\zeroat{B}} &
  {B} }$ depicts the unique morphisms to and from \zeroo{} and
$\zeroa{A}{B} = \comp{\zeroat{B}}{\zeroaf{A}}$.

\msection{Subobjects in extensive contexts}
\label{sec:subobj-extens-categ}

Let $\mathcal{A} = (\Bb{A}, \mathsf{E}, \mathsf{M})$ be an extensive
context, i.e., \Bb{A}{} is \emph{lextensive}
\cite[see][]{CarboniLackWalters1993}{.}  Extensivity ensures the
initial object is strict \cite[see][\S2]{CarboniLackWalters1993},
while the \fact{\textsf{E}}{\textsf{M}} ensures coproduct injections
are admissible monomorphisms. To see this, since the coproduct
$\xymatrix{ {\termo} \ar[r]^-{\iota_1} & {\termo+\termo}
  \ar@{<-}[r]^-{\iota_2} & {\termo} }$ exists, $\iota_i$ (for
$i = 1, 2$) is a split monomorphism and hence is in \textsf{M}. For
any coproduct
$\xymatrix{ {X} \ar[r]^-{\iota_X} & {X+Y} \ar@{<-}[r]^-{\iota_Y} & {Y} }$,
using extensivity, since
$\iota_X = \finv{(\terma{X}+\terma{Y})}{\iota_1}$ and
$\iota_Y = \finv{(\terma{X}+\terma{Y})}{\iota_2}$ the coproduct injections are
admissible monomorphisms. Consequently from $\termo = \inito+\termo$,
every extensive context is admissibly quasi-pointed (i.e., a
quasi-pointed category in which \Arr{\inita{\termo}}{\inito}{\termo}{}
is an admissible monomorphism, \cite[see][for
quasi-pointed]{Bourn2001, GoswamiJanelidze2017} and \cite[see][Remark
(O)]{2021-clos} with a strict initial object. Hence equivalently, every
morphism reflects zero, or equivalently every preneighbourhood
morphism is continuous \cite[see][\S9]{2021-clos}. Furthermore, in a
lextensive category
\begin{equation*}
  \xymatrixcolsep{9em}
  \xymatrix{
    {\opp{\Bb{A}} \times \opp{\Bb{A}}} \ar@<1em>[r]^-{\Sub{-+-}{}}="u"
    \ar@<-1em>[r]_-{\Sub{-}{}\times\Sub{-}{}}="d" & {\PosZero}
    \ar <0pt,-0.3em>+"u";<0pt,0.3em>+"d"^-{\iota}
  }
\end{equation*}
is a natural isomorphism, where
$\iota_{X,Y} = \opair{\finv{\iota_X}{}}{\finv{\iota_Y}{}}$ and
$\inv{\iota_{X,Y}} = +_{X,Y}$ the sum of subobjects and \Sub{X}{}{} is the
partially ordered set of all subobjects of $X$.

\begin{Lemma}
  \label{lem:extensivity-sub-of-sum}
  In a lextensive category \Bb{A}{,} let
  $\xymatrix{ {\Sub{X}{}} \ar[r]^-{L_X} & {\Sub{X+Y}{}}
    \ar@{<-}[r]^-{R_Y} & {\Sub{Y}{}} }$ ($X, Y \in \Bb{A}_0$) define the
  functions
  \begin{align}
    \label{eq:left-addition}
    L_{X}(m) & = m + \sigma_{Y}, & \text{ for }m \in\Sub{X}{}, \\[2pt]
    \label{eq:right-addition}
    R_{Y}(n) & = \sigma_{X} + n, & \text{ for }n \in\Sub{Y}{}.
  \end{align}
  Then the morphisms in the diagram:
  \begin{equation}
    \label{eq:nearly-biproduct}
    \xymatrix{
      {\Sub{X}{}} \ar@<1ex>[r]^-{L_X} \ar@{<-}@<-1ex>[r]_-{\finv{\iota_X}{}} &
      {\Sub{X+Y}{}} \ar@<1ex>@{<-}[r]^-{\iota_Y}
      \ar@<-1ex>[r]_-{\finv{\iota_Y}{}} &
      {\Sub{Y}{}}
    }
  \end{equation}
  in \PosZero{} satisfy the equations:
  \begin{align}
    \label{eq:left-rules}
    \comp{\finv{\iota_X}{}}{L_X} & = \id{\Sub{X}{}},
    & & \comp{\finv{\iota_X}{}}{R_Y} & = \zeroa{\Sub{Y}{}}{\Sub{X}{}},
    \\[2pt]
    \comp{\finv{\iota_Y}{}}{L_Y} & = \id{\Sub{Y}{}},
    & & \comp{\finv{\iota_Y}{}}{L_X} & = \zeroa{\Sub{X}{}}{\Sub{Y}{}},
  \end{align}
  and
  \begin{equation}
    \label{eq:sum-rule}
    (\comp{L_X}{\finv{\iota_X}{}}) \vee (\comp{R_Y}{\finv{\iota_Y}{}}) =
    \id{\Sub{X+Y}{}}.
  \end{equation}
  In particular, \Sub{X+Y}{} is a biproduct of \Sub{X}{} and
  \Sub{Y}{}.
\end{Lemma}

Recall: since \Arr{U}{\pNHD{\Bb{A}}}{\Bb{A}}{} is topological
\cite[see][Theorem 4.8(a)]{2020}, the coproduct $X + Y$ of internal
preneighbourhood spaces \opair{X}{\mu}{,} \opair{Y}{\phi}{} is given the
largest preneighbourhood system $\mu+\phi$ such that the coproduct
injections $\iota_X$, $\iota_Y$ are preneighbourhood morphisms.

\begin{Thm}
  \label{thm:joinleftadjttosuminjpb}
  Let
  $\xymatrix{ {\opair{X}{\mu}} \ar@{>->}[]!<1.92em,0em>;[r]^-{\iota_X}
    & {\opair{X+Y}{\mu+\phi}} & {\opair{Y}{\phi}}
    \ar@{>->}[]!<-1.8em,0em>;[l]_-{\iota_Y} }$ be the coproduct of the
  internal preneighbourhood spaces \opair{X}{\mu} and \opair{Y}{\phi}.
  The isomorphism $\iota$ restrict to adjunctions:
  \begin{equation}
    \label{eq:adm-subobj-adjt}
    \xymatrixcolsep{10.8em}
    \xymatrix{
      {\Sub{X}{\mathsf{M}}\times\Sub{Y}{\mathsf{M}}}
      \ar@<1.44ex>[r]^-{\img{\iota_X}{}\vee\img{\iota_Y}{}}
      \ar@<-1.44ex>@{<-}[r]_-{\iota_{X,Y}=\opair{\finv{\iota_X}{}}{\finv{\iota_Y}{}}}
      \ar@{}[r]|-{\bot}
      & {\Sub{X+Y}{\mathsf{M}}}
    },
  \end{equation}
  and
  \begin{equation}
    \label{eq:closed-subobj-adjt}
    \xymatrixcolsep{10.8em}
    \xymatrix{
      {\mc{\mu}\times\mc{\phi}}
      \ar@<1.44ex>[r]^-{\cls{(\img{\iota_X}{}\vee\img{\iota_Y}{})}{\mu+\phi}}
      \ar@<-1.44ex>@{<-}[r]_-{\iota_{X,Y}=\opair{\finv{\iota_X}{}}{\finv{\iota_Y}{}}}
      \ar@{}[r]|-{\bot}
      & {\mc{\mu+\phi}}
    }
  \end{equation}
  in \PosZero{}.
\end{Thm}

\begin{proof}
  Given any $m \in \Sub{X}{\mathsf{M}}$, $n \in \Sub{Y}{\mathsf{M}}$
  and $p \in \Sub{X + Y}{\mathsf{M}}$:
  \begin{equation*}
    \begin{aligned}
      (m, n) \leq \iota_{X, Y}(p) & \Leftrightarrow m \leq
      \finv{\iota_X}{p}\text{ and }n \leq \finv{\iota_Y}{p} \\
      & \Leftrightarrow \img{\iota_X}{m} \leq p
      \text{ and }\img{\iota_Y}{n} \leq p \\
      & \Leftrightarrow \img{\iota_X}{m} \vee \img{\iota_Y}{n} \leq p,
    \end{aligned}
  \end{equation*}
  completing the proof the adjunction \adjt{(\img{\iota_X}{} \vee
    \img{\iota_Y}{})}{\iota_{X, Y}}.  Given $u \in \mc{\mu}$,
  $v \in \mc{\phi}$ and $w \in \mc{\mu + \phi}$:
  \begin{equation*}
    \begin{aligned}
      (u, v) \leq \iota_{X, Y}(w) & \Leftrightarrow
      u \leq \finv{\iota_X}{w}\text{ and }v \leq \finv{\iota_Y}{w} \\
      & \Leftrightarrow \img{\iota_X}{u} \leq w \text{ and }
      \img{\iota_Y}{v} \leq w \\
      & \Leftrightarrow \img{\iota_X}{u} \vee \img{\iota_Y}{v} \leq w \\
      & \Leftrightarrow \cls{(\img{\iota_X}{u} \vee \img{\iota_Y}{v})}
      {\mu + \phi} \leq w,
    \end{aligned}
  \end{equation*}
  completing the proof of the adjunction \adjt{\cls{(\img{\iota_X}{}
      \vee \img{\iota_Y}{})}{\mu + \phi}}{\iota_{X, Y}}.

\end{proof}

\msection{\fact{\textsf{E}}{\textsf{M}} of finite sums}
\label{sec:fact-sum-morphisms}

Given
$\xymatrix{ {X} \ar[r]^-{\iota_X} & {X + Y} \ar@{<-}[r]^-{\iota_Y} & {Y} }$,
the coproduct of $X$ and $Y$, morphisms \Arr{f}{X}{Z} and
\Arr{g}{Y}{Z}, the unique morphism from $X+Y$ to $Z$ is
\Arr{\Opair{f}{g}}{X+Y}{Z}, where $f = \comp{\Opair{f}{g}}{\iota_X}$ and
$g = \comp{\Opair{f}{g}}{\iota_Y}$.

\begin{Thm}
  \label{sumisadmissible-equivalents}
  In an extensive context, given morphisms \Arr{a}{A}{X} and
  \Arr{b}{B}{Y} if the \fact{\textsf{E}}{\textsf{M}} of
  \Opair{\img{\iota_X}{a}}{\img{\iota_Y}{b}}{} is
  $\Opair{\img{\iota_X}{a}}{\img{\iota_Y}{b}} = \comp{\bigl(\img{\iota_X}{a} \vee
    \img{\iota_Y}{b}\bigr)}{e}$ then
  \begin{equation}
    \label{emfactforsum-eq}
    a + b = \comp{\bigl(\img{\iota_X}{a} \vee
      \img{\iota_Y}{b}\bigr)}{\comp{e}{\bigl(\rest{\iota_X}{A} +
        \rest{\iota_Y}{B}\bigr)}},
  \end{equation}
  is the \fact{\mathsf{E}}{\mathsf{M}} of $a + b$.  In particular,
  $a + b \in \mathsf{M}$ if and only if
  $a + b = \img{\iota_X}{a} \vee \img{\iota_Y}{b}$.
\end{Thm}

\begin{proof}
  \begin{figure}
    \centering
    \begin{equation*}
      % \resizebox{36em}{!}{
      \xymatrixcolsep{2.4em}
      \xymatrixrowsep{3em}
      \xymatrix{
        {A} \ar[ddd]_{a} \ar@{>->}[]!<1.2em,0em>;[rrrr]^-{\iota_A}
        \ar@{->>}[dr]|{\rest{\iota_X}{A}} & & & &
        {A + B} \ar[ddd]|-{a + b} 
        \ar@{->>}[dll]|
        {\rest{\iota_X}{A} + \rest{\iota_Y}{B}} & &
        {B} \ar[ddd]^-{b} \ar@{>->}[]!<-1.2em,0em>;[ll]_-{\iota_B}
        \ar@{->>}[dl]|{\rest{\iota_Y}{B}} \\
        & {\img{\iota_X}{A}} \ar[r]^-{u_X}
        \ar@/_2.4ex/@{>->}[]!<12pt,-12pt>;[ddrrr]|{\img{\iota_X}{a}}
        & {\img{\iota_X}{A} + \img{\iota_Y}{B}}
        \ar@{<-}[rrr]|(0.36){u_Y}|(0.738){\hole}
        \ar@{->>}[dr]|{e} & & & {\img{\iota_Y}{B}}
        \ar@/^{2ex}/@{>->}[]!<-9pt,-12pt>;[ddl]|
        {\img{\iota_X}{b}} & \\
        & & & {\img{\iota_X}{A} \vee \img{\iota_Y}{B}}
        \ar@{>->}[]!<12pt,-12pt>;[dr]|{\img{\iota_X}{a} \vee
          \img{\iota_Y}{b}} & & & \\
        {X} \ar@{>->}[]!<1.2em,0em>;[rrrr]_-{\iota_X} & & & &
        {X + Y}  &
        &
        {Y} \ar@{>->}[]!<-1.2em,0em>;[ll]^-{\iota_Y} \\    
      }
      % }
    \end{equation*}    
    \caption{Factorising sum of morphisms}
    \label{fig:sumisadmissible-equivalents}
  \end{figure}
  Consider the diagram in Figure
  \ref{fig:sumisadmissible-equivalents}{,} where the top and the
  bottom rows are coproducts as shown. The second row from the top
  exhibits the coproduct of \img{\iota_X}{A} and \img{\iota_Y}{B} and
  the composition $\comp{(\img{\iota_X}{a} \vee \img{\iota_Y}{b})}{e}$
  provides the \fact{\mathsf{E}}{\mathsf{M}} of
  \Opair{\img{\iota_X}{a}}{\img{\iota_Y}{b}}{}. Furthermore, since the
  two rows from the top are both coproducts there exists the unique
  morphism $\rest{\iota_X}{A} + \rest{\iota_Y}{B}$ from $\mathsf{E}$
  which makes the top squares to commute.  Hence
  $a + b = \comp{\bigl(\img{\iota_X}{a} \vee
    \img{\iota_Y}{b}\bigr)}{\comp{e}{\bigl(\rest{\iota_X}{A} +
      \rest{\iota_Y}{B}\bigr)}}$ is the \fact{\mathsf{E}}{\mathsf{M}}
  of $a + b$.  In particular, if $a + b \in \mathsf{M}$ then
  $\comp{e}{\bigl(\rest{\iota_X}{A} + \rest{\iota_Y}{B})} \in
  \mathsf{E} \cap \mathsf{M} = \Iso{\Bb{A}}$ and hence
  $a + b = \img{\iota_X}{a} \vee \img{\iota_Y}{b}$, proving the
  \ul{\em only if} part of the second statement. The \ul{\em if} part
  of the statement is obviously trivial.
\end{proof}

As an immediate application of this along with Theorem
\ref{thm:joinleftadjttosuminjpb}:
\begin{Cor}
  \label{sumasinverse}
  The order preserving maps \Arr{\iota_{X, Y}}{\Sub{X +
      Y}{\mathsf{M}}}{\Sub{X}{\mathsf{M}} \times \Sub{Y}{\mathsf{M}}}
  is an isomorphism of partially ordered sets if and only if
  $+_{X, Y} = \img{\iota_X}{} \vee \img{\iota_Y}{} = \inv{\iota_{X,
      Y}}$.
\end{Cor}

\msection{Sum of admissible subobjects}
\label{sec:enable-iso-of-sums}

Recall from \cite[][\S2{}]{Janel1997b}, given morphisms $e$ and $m$,
\down{e}{m}{} if $\comp{v}{e} = \comp{m}{u}$ implies the existence of
a unique morphism $w$ such that $v = \comp{m}{w}$ and
$u = \comp{w}{e}$. In the \fact{\textsf{E}}{\textsf{M}},
$\mathsf{E} = \mathsf{M}^{\uparrow} = \bigl\{ x: m \in \mathsf{M} \Rightarrow \down{x}{m}
\bigr\}$ and
$\mathsf{M} = \mathsf{E}^{\downarrow} = \bigl\{ x: e \in \mathsf{E} \Rightarrow \down{e}{x}
\bigr\}$.

\begin{Thm}
  \label{Extensivity>pbstabmonicformalepis<>admsumclosed}
  In an extensive context, finite sum of admissible subobjects is an
  admissible subobject if and only if the monomorphisms in
  $\mathsf{E}$ between finite sums are stable under pullbacks along
  coproduct injections.
\end{Thm}

\begin{proof}
  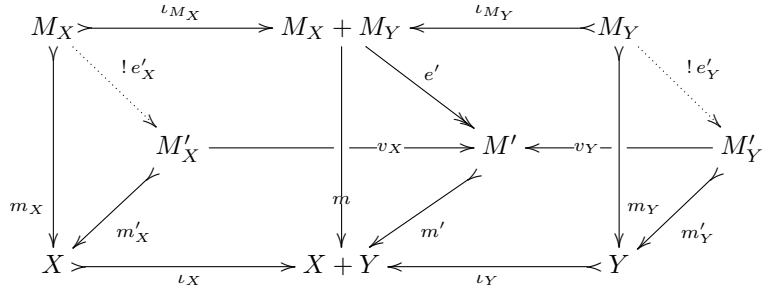
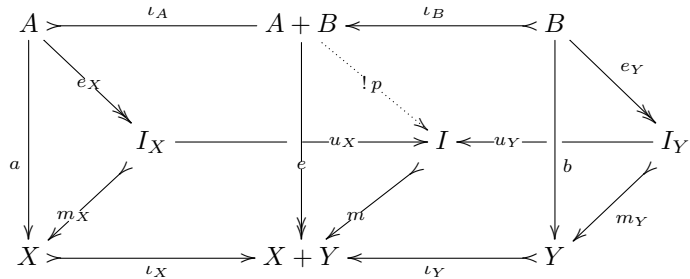
\begin{figure}[b]
    \centering
    \begin{subfigure}{0.96\linewidth}
      \begin{equation*}
        % \resizebox{0.9\textwidth}{0.24cm}{
        \xymatrixcolsep{2.4em}
        \xymatrixrowsep{3em}
        \xymatrix{
          {M_X} \ar@{>->}[]!<0em,-1.2em>;[dd]_(0.72){m_X}
          \ar@{>->}[]!<1.44em,0em>;[rr]^-{\iota_{M_X}}
          \ar@{.>}[dr]^-{!\, e'_X} & & {M_{X}+M_{Y}}
          \ar[dd]|(0.72){m} 
          \ar@{->>}[dr]^-{e'} &
          & {M_Y} \ar@{>->}[]!<0em,-1.2em>;[dd]^(0.72){m_Y}
          \ar@{.>}[dr]^-{!\, e'_Y}
          \ar@{>->}[]!<-1.44em,0em>;[ll]_-{\iota_{M_Y}} & \\
          & {M'_X} \ar@{>->}[]!<-1.2em,-1.2em>;[dl]^-{ m'_X}
          \ar[rr]|-{\hole}|(0.66){v_X} & & {M'}
          \ar@{>->}[]!<-1.2em,-1.2em>;[dl]^-{m'}
          \ar@{<-}[rr]|(0.36){v_Y}|-{\hole} & & {M'_Y}
          \ar@{>->}[]!<-1.2em,-1.2em>;[dl]^-{m'_Y} \\
          {X} \ar@{>->}[]!<1.2em,0em>;[rr]_-{\iota_X} & & {X + Y} 
          & & {Y} \ar@{>->}[]!<-1.2em,0em>;[ll]^-{\iota_Y}&
        }
        % }
      \end{equation*}
      \caption{\emph{if} part}
      \label{fig:sum-adm-sub-if-part}
    \end{subfigure}

    \begin{subfigure}{0.96\linewidth}
      \begin{equation*}
        % \resizebox{0.9\textwidth}{0.24cm}{
        \xymatrixcolsep{3em}
        \xymatrixrowsep{3em}
        \xymatrix{
          {A} \ar@{>-}[]!<1.2em,0em>;[rr]^-{\iota_A}
          \ar[dd]_(0.6){a}
          \ar@{->>}[dr]|{e_X} & &
          {A + B} 
          \ar@{->>}[dd]|(0.6){e}
          \ar@{.>}[dr]|{!\, p}
          & &
          {B} \ar[dd]^(0.6){b} \ar@{>->}[]!<-1.2em,0em>;[ll]_-{\iota_B}
          \ar@{->>}[dr]^-{e_Y} &
          \\
          & {I_X} \ar@{>->}[]!<-1.2em,-1.2em>;[dl]|{m_X}
          \ar[rr]|(0.492){\hole}|(0.66){u_X} & &
          {I}
          \ar@{<-}[rr]|(0.282){u_Y}|(0.48){\hole}
          \ar@{>->}[]!<-1.2em,-1.2em>;[dl]|{m}
          & &
          {I_Y} \ar@{>->}[]!<-1.2em,-1.2em>;[dl]^-{m_Y} \\
          {X} \ar@{>->}[]!<1.2em,0em>;[rr]_-{\iota_X} & & {X + Y}
          & & {Y} \ar@{>->}[]!<-1.2em,0em>;[ll]^-{\iota_Y}& \\    
        }
        % }
      \end{equation*}
      \caption{\emph{only if} part}
      \label{fig:sum-adm-sub-onlyif-part}
    \end{subfigure}
    
    \caption{Sum of admissible subobjects}
    \label{fig:sum-adm-sub}
  \end{figure}
  Towards proof of the \emph{if} part, given the admissible subobjects
  $m_X \in \Sub{X}{\mathsf{M}}$ and $m_Y \in \Sub{Y}{\mathsf{M}}$ consider
  the diagram in Figure \ref{fig:sum-adm-sub-if-part} in which both
  the rows are coproduct diagrams as shown, $m = m_X + m_Y$. Hence the
  vertical squares are all pullbacks from extensivity.  Since sum of
  monomorphisms is a monomorphism, $m \in \Mono{\Bb{A}}$.  Let
  $m = \comp{m'}{e'}$ be the \fact{\mathsf{E}}{\mathsf{M}} of
  $m$. Hence $e' \in \mathsf{E} \cap \Mono{\Bb{A}}$.  The morphisms
  $m'_X$ and $m'_Y$ are the pullbacks $m'$ along $\iota_X$ and
  $\iota_Y$ respectively. Using extensivity, $m'_X + m'_Y = m'$. Hence
  $e'$ is a morphism between finite sums. Using property of pullback
  squares there exist unique morphisms $e'_X$ and $e'_Y$ such that
  $\comp{m'_X}{e'_X} = m_X$, $\comp{m'_Y}{e'_Y} = m_Y$,
  $\comp{v_X}{e'_X} = \comp{e'}{\iota_{M_X}}$ and
  $\comp{v_Y}{e'_Y} = \comp{e'}{\iota_{M_Y}}$.  Consequently
  $e'_X, e'_Y \in \mathsf{M}$ and all the squares of the diagram are
  pullback squares.  From hypothesis,
  $e'_X, e'_Y \in \mathsf{E} \cap \mathsf{M} = \Iso{\Bb{A}}$.  Hence
  $M = M_X + M_Y \approx M_{X'} + M_{Y'} = M'$ forcing
  $e' \in \Iso{\Bb{A}}$, proving \ul{\em if} part of the theorem.
  Towards proof of the \emph{only if} part, the hypothesis along with
  Corollary \ref{sumasinverse}{} ensures for each $X, Y \in \Bb{A}_0$,
  $\inv{\iota_{X, Y}} = \img{\iota_X}{} \vee \img{\iota_Y} = +_{X, Y}$.  Consider the
  diagram in Figure \ref{fig:sum-adm-sub-onlyif-part}{} where the top
  and bottom row are coproduct diagrams as shown and
  $e \in \mathsf{E} \cap \Mono{\Bb{A}}$ is a morphism between finite sums.
  Since $e$ is a monomorphism, there exist unique monomorphisms
  $a \in \Sub{X}{}$ and $b \in \Sub{Y}{}$ such that $e = a + b$; using
  extensivity $\finv{\iota_X}{e} = a$ and $\finv{\iota_Y}{e} = b$.  Let
  $a = \comp{m_X}{e_X}$ and $b = \comp{m_Y}{e_Y}$ be the
  \fact{\mathsf{E}}{\mathsf{M}} for $a$ and $b$ respectively. From
  hypothesis, $m = m_X + m_Y \in \Sub{X+Y}{\mathsf{M}}$. Hence using
  extensivity, the horizontal squares are both pullback squares,
  $I = I_X + I_Y$, $u_X$ and $u_Y$ are coproduct injections.  Since
  the top row is a coproduct, there exists the unique morphism $p$
  such that $\comp{p}{\iota_A} = \comp{u_X}{e_X}$,
  $\comp{p}{\iota_B} = \comp{u_Y}{e_Y}$ and $p = e_X + e_Y$.  Hence,
  the slanting top squares are pullback squares.  Thus, from the front
  vertical squares $e = a + b = \comp{m}{p}$.  Since
  $e \in \mathsf{E}$ and \down{e}{m}, $m$ is an isomorphism. Hence,
  $m_X = \finv{\iota_X}{m}$ and $m_Y = \finv{\iota_Y}{m}$ are also
  isomorphisms, proving $a, b \in \mathsf{E}$, completing the proof.
\end{proof}

\msection{Sum of closed embeddings}
\label{sec:sum-of-closed-morphisms}

Every lextensive category is distributive \cite[see][Proposition
4.5]{CarboniLackWalters1993}. Consequently, each lattice
\Sub{X}{\mathsf{M}}{} of admissible subobjects is distributive. Hence
the closure operator
\Arr{\cls{}{\mu}}{\Sub{X}{\mathsf{M}}}{\Sub{X}{\mathsf{M}}}{} on each
internal preneighbourhood space \opair{X}{\mu}{} is additive
\cite[see][Theorem 3.1(d), Remark (C)]{2020}. Recall: any
preneighbourhood morphism \Arr{f}{\opair{X}{\mu}}{\opair{Y}{\phi}} factors
as
$\xymatrixcolsep{3em}\xymatrix{ {X} \ar@{->>}[r]^-{f^{\mathsf{E}}} &
  {\Img{f}} \ar@{>->}[]!<2.16em,0em>;[r]^-{u_f}
  \ar@/_{4.2ex}/@{>->}[]!<0em,-0.96em>;[rr]|{f^{\mathsf{M}}} &
  {\overline{\Img{f}}}
  \ar@{>->}[]!<2.16em,0em>;[r]^-{\cls{f^{\mathsf{M}}}{\phi}} & {Y} }$,
where $f = \comp{f^{\mathsf{M}}}{f^{\mathsf{E}}}$ is the
\fact{\textsf{E}}{\textsf{M}} of $f$, \comp{u_f}{f^{\mathsf{E}}}{} is
a dense morphism and
$f = \comp{\cls{f^{\mathsf{M}}}{\phi}}{(\comp{u_f}{f^{\mathsf{E}}})}$ is
the dense-(closed embedding) factorisation of $f$ \cite[see][proof of
Theorem 5.1]{2021-clos}.

\begin{Thm}
  \label{sumofclosedsubobjects}
  In any extensive context the following are equivalent:
  \begin{enumerate}[label=(\alph*)]
  \item \label{item:clemb-fin-sum-closed} Every finite sum of closed
    embeddings is a closed embedding.

  \item \label{item:clemb-fin-sum-is-adm-coprdinj-closed} Every finite
    sum of admissible subobjects is an admissible subobject and each
    coproduct injection is a closed embedding.

  \item \label{item:dense-sb-stable-along-coprdinj}{} Each dense
    morphism between finite sums is stable under pullbacks along
    coproduct injections.

  \end{enumerate}

\end{Thm}

\begin{proof}
  \begin{figure}
    \centering
    \begin{subfigure}{0.96\linewidth}
      \begin{equation*}
        \xymatrixcolsep{3em}
        \xymatrixrowsep{3em}
        \xymatrix{
          {Z_1} \ar[rr]^-{i_{1}} \ar[dd]_-{d_{1}} \ar[dr]^-{h_1}
          & & {Z} \ar@{<-}[rr]^-{i_{2}} \ar[dd]|(0.6){d}
          \ar@{.>}[dr]|{!\,\Opair{h_1}{h_2}}
          & & {Z_{2}} \ar[dd]|(0.6){d_{2}} \ar[dr]^-{h_2} \\
          & {M_{1}} \ar@{>->}[]!<1.2em,0em>;[rr]|(0.396){\hole}|(0.6){\iota_1}
          \ar@{>->}[]!<-0.84em,-0.84em>;[dl]^-{m_1}
          & & {M_{1} + M_{2}} 
          \ar@{>->}[]!<-0.84em,-0.96em>;[dl]!<1.8em,0.84em>|{m_1+m_2}
          & & {M_{2}} \ar@{>->}[]!<-0.84em,-0.84em>;[dl]^-{m_2}
          \ar@{>->}[]!<-1.2em,0em>;[ll]|(0.36){\hole}|(0.54){\iota_2} \\
          {X} \ar@{>->}[]!<1.2em,0em>;[rr]_-{\iota_X} & & {X + Y}
          & & {Y} \ar@{>->}[]!<-1.2em,0em>;[ll]^-{\iota_Y}
        }
      \end{equation*}
      \caption{\ref{item:clemb-fin-sum-closed}{} implies
        \ref{item:dense-sb-stable-along-coprdinj}}
      \label{fig:a2c}{}
    \end{subfigure}

    \begin{subfigure}{0.96\linewidth}
      \begin{equation*}
        \xymatrixcolsep{3em}
        \xymatrixrowsep{3em}
        \xymatrix{
          {A} \ar@{>->}[]!<1.2em,0em>;[rr]^-{\iota_A}
          \ar@{>->}[]!<0em,-0.96em>;[dd]_-a
          \ar@{.>}[dr]^-{!\,u}
          & & {A+B}  \ar[dd]|(0.6){a+b}
          \ar[dr]|{w}
          & & {B} \ar@{>->}[]!<0em,-0.96em>;[dd]|(0.6){b}
          \ar@{.>}[dr]^-{!\,v} \ar@{>->}[]!<-1.2em,0em>;[ll]_-{\iota_B}\\
          & {M} \ar[rr]|(0.432){\hole}|(0.6){\iota_1}
          \ar@{>->}[]!<-0.84em,-0.84em>;[dl]^-{a'}
          & & {\overline{\Img{a+b}}}
          \ar@{<-}[rr]|(0.36){\iota_2}|(0.564){\hole}
          \ar@{>->}[]!<-2.4em,-1.2em>;[dl]!<1.8em,0.84em>
          |(0.42){\cls{(a+b)^{\mathsf{M}}}{\mu+\phi}}
          & & {N} \ar@{>->}[]!<-0.84em,-0.84em>;[dl]^-{b'} \\
          {X} \ar@{>->}[]!<1.2em,0em>;[rr]_-{\iota_X} & & {X + Y}
          & & {Y} \ar@{>->}[]!<-1.2em,0em>;[ll]^-{\iota_Y} 
        }
      \end{equation*}
      \caption{\ref{item:dense-sb-stable-along-coprdinj}{} implies
        \ref{item:clemb-fin-sum-closed}{}}
      \label{fig:c2a}{}
    \end{subfigure}
    \caption{Equivalence of \ref{item:dense-sb-stable-along-coprdinj}
      and \ref{item:clemb-fin-sum-closed}}
    \label{fig:sum-closed-emb}
  \end{figure}
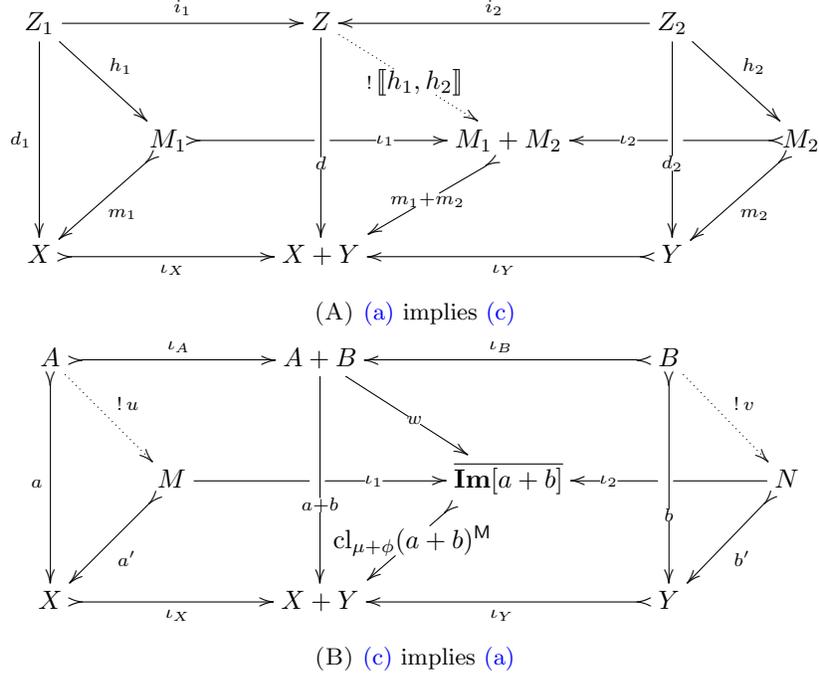
  
  Let
  $\xymatrix{ {\opair{X}{\mu}} \ar[r]^-{\iota_X} & {\opair{X+Y}{\mu+\phi}}
    \ar@{<-}[r]^-{\iota_Y} & {\opair{Y}{\phi}} }$ be the coproduct of
  internal preneighbourhood spaces \opair{X}{\mu}{} and \opair{Y}{\phi}{.}
  Assuming \ref{item:clemb-fin-sum-closed}{,} for each
  $a \in \mc{\mu}$ and $b \in \mc{\phi}$ there exists a unique
  $c \in \mc{\mu + \phi}$ such that $a = \finv{\iota_X}{c}$ and
  $b = \finv{\iota_Y}{c}$.  Using Theorem
  \ref{sumisadmissible-equivalents} and extensivity,
  $c = a + b = \img{\iota_X}{a} \vee \img{\iota_Y}{b}$. Hence from Theorem
  \ref{thm:joinleftadjttosuminjpb},
  $\inv{\iota_{X, Y}} = \cls{(- + -)}{\mu + \phi}$.  In particular, for any
  $a \in \mc{\mu}$,
  $\cls{\img{\iota_X}{a}}{\mu + \phi} = \cls{(a + \sigma_{Y})}{\mu + \phi} = a + \sigma_{Y} =
  \img{\iota_X}{a}$, proving $\iota_X$ is a closed morphism; similarly,
  $\iota_Y$ is a closed morphism. This proves
  \ref{item:clemb-fin-sum-is-adm-coprdinj-closed}{.} Assuming
  \ref{item:clemb-fin-sum-is-adm-coprdinj-closed}{,} if the coproduct
  injections $\iota_X$ and $\iota_Y$ are closed embeddings in \Bb{A} and
  $a \in \mc{\mu}, b \in \mc{\phi} \Rightarrow a + b \in \mathsf{M}$ then for each
  $(a, b) \in \mc{\mu} \times \mc{\phi}$, using Theorem
  \ref{sumisadmissible-equivalents} and additivity of the closure
  operator,
  $\cls{(a + b)}{\mu + \phi} = \cls{(\img{\iota_X}{a} \vee \img{\iota_Y}{b})}{\mu + \phi} =
  (\cls{\img{\iota_X}{a}}{\mu+\phi}) \vee (\cls{\img{\iota_Y}{b}}{\mu+\phi}) =
  (\img{\iota_X}{\cls{a}{\mu}}) \vee (\img{\iota_Y}{\cls{b}{\phi}}) = \img{\iota_X}{a} \vee
  \img{\iota_Y}{b} = a + b$, proving \ref{item:clemb-fin-sum-closed}{.}
  
  Assuming \ref{item:clemb-fin-sum-closed}{,} let
  \Arr{d}{\opair{Z}{\psi}}{\opair{X+Y}{\mu+\phi}}{} be a dense
  morphism. Consider the diagram in Figure \ref{fig:a2c}{} in which
  $d$ is pulled back along the coproduct injections $\iota_X$ and
  $\iota_Y$ producing $d_{1}$ and $d_{2}$. Using extensivity,
  $Z = Z_{1} + Z_{2}$ with coproduct injections $i_{1}$, $i_{2}$ and
  $d = d_{1} + d_{2}$. If $d_{i} = \comp{m_i}{h_i}$, where $m_{i}$ are
  closed embeddings ($i = 1, 2$) then since $Z = Z_{1} + Z_{2}$, there
  exists the unique morphism \Opair{h_1}{h_2}{} making the whole
  diagram to commute. From assumption $m_{1} + m_{2}$ is a closed
  embedding. Since $d = \comp{(m_1+m_2)}{\Opair{h_1}{h_2}}$, and $d$
  is dense, $m_{1} + m_{2}$ is an isomorphism. Hence
  $m_{1} = \finv{\iota_X}{(m_1+m_2)}$ and
  $m_{2} = \finv{\iota_Y}{(m_1+m_2)}$ are both isomorphisms, proving
  $d_{1}$ and $d_{2}$ are dense. Conversely, assuming
  \ref{item:dense-sb-stable-along-coprdinj}{,} given $a \in \mc{\mu}$,
  $b \in \mc{\phi}$, consider the diagram in Figure \ref{fig:c2a}{} where
  the top and bottom row of the front vertical squares are both
  coproducts. Let
  $a + b = \comp{\cls{(a+b)^{\mathsf{M}}}{\mu+\phi}}{w}$ be the
  dense-(closed embedding) factorisation of $a + b$. The closed
  embedding \cls{(a+b)^{\textsf{M}}}{\mu+\phi}{} is pulled back along the
  coproduct injections to obtain the morphisms $a'$ and $b'$. Hence,
  from extensivity, $\cls{(a+b)^{\mathsf{M}}}{\mu+\phi} = a' + b'$,
  $\overline{\Img{a+b}} = M + N$ with coproduct injections $\iota_1$ and
  $\iota_2$. Furthermore, since preneighbourhood morphisms are continuous,
  $a' = \finv{\iota_X}{\cls{(a+b)^{\mathsf{M}}}{\mu+\phi}}$,
  $b' = \finv{\iota_Y}{\cls{(a+b)^{\mathsf{M}}}{\mu+\phi}}$ are both closed
  embeddings. Since
  $\comp{\cls{(a+b)^{\mathsf{M}}}{\mu+\phi}}{\comp{w}{\iota_A}} =
  \comp{\iota_X}{a}$ and
  $\comp{\cls{(a+b)^{\mathsf{M}}}{\mu+\phi}}{\comp{w}{\iota_B}} =
  \comp{\iota_Y}{b}$, from the horizontal pullback squares there exist
  unique morphisms $u$, $v$ such that the whole diagram
  commutes. Hence, from extensivity, all the squares in the diagram
  are pullback squares. From hypothesis, $u$ and $v$ are dense
  morphisms. Since $a'$ is closed, and closure operation is transitive
  \cite[see][Theorem 3.1(e)]{2021-clos},
  $\comp{a'}{u} = a = \cls{a}{\mu} =
  \comp{a'}{\cls{u}{\rest{\mu}{M}}}$, implying $u$ is a closed
  embedding too. Similarly, $v$ is a closed embedding. Since dense
  closed embeddings are isomorphisms \cite[see][Theorem
  5.1(b)]{2021-clos}, $u$ and $v$ are isomorphisms. Consequently,
  $w = u + v$ is an isomorphism, proving $a + b$ to be a closed
  embedding.
\end{proof}

\begin{Cor}
  \label{cor:sum-of-closed-morphisms}
  In an extensive context with finite sum of admissible subobjects an
  admissible subobject, finite sum of closed morphisms is a closed
  morphism if and only if, the coproduct injections are closed.
\end{Cor}

\begin{proof}
  The \emph{only if} part follows from Theorem and description of
  closed embeddings \cite[][Theorem 4.1(b), second
  part]{2021-clos}. For the \emph{if} part, if
  \Arr{f}{\opair{A}{\alpha}}{\opair{X}{\mu}}{} and
  \Arr{g}{\opair{B}{\beta}}{\opair{Y}{\phi}}{} are closed morphisms,
  $p = p_{A} + p_{B} \in \mc{\alpha+\beta}$, where
  $p_{A} = \finv{\iota_A}{p}$ and $p_{B} = \finv{\iota_B}{p}$, then using
  equation \eqref{eq:clos-sum-is-sum-of-clos} in Lemma below and
  Corollary \ref{factsumcommute2}{} in next section,
  $\img{f+g}{\cls{p}{\alpha+\beta}} = \img{f+g}{\cls{(p_A+p_B)}{\alpha+\beta}} =
  \img{f+g}{(\cls{p_A}{\alpha}+\cls{p_B}{\beta})} = \img{f}{\cls{p_A}{\alpha}} +
  \img{g}{\cls{p_B}{\beta}} = \cls{\img{f}{p_A}}{\mu} +
  \cls{\img{f}{p_B}}{\phi} = \cls{(\img{f}{p_A}+\img{g}{p_B})}{\mu+\phi} =
  \cls{\img{f+g}{p}}{\mu+\phi}$, proving $f + g$ is a closed morphism using
  \cite[][Theorem 4.1(b)]{2021-clos}.
\end{proof}

\begin{Lemma}
  \label{lem:clos-sum-is-sum-of-clos}
  If finite sum of closed embeddings is a closed embedding then for
  every internal preneighbourhood spaces \opair{X}{\mu}{,}
  \opair{Y}{\phi}{,} $a \in \mc{\mu}$ and $b \in \mc{\phi}$:
  \begin{equation}
    \label{eq:clos-sum-is-sum-of-clos}
    \cls{(a+b)}{\mu+\phi} = \cls{a}{\mu} + \cls{b}{\phi}.
  \end{equation}
\end{Lemma}

\begin{proof}
  Since from Theorem the coproduct injections are closed, finite sum
  of admissible subobjects are admissible, using Theorem
  \ref{thm:joinleftadjttosuminjpb}:
  \begin{multline*}
    \cls{a}{\mu} + \cls{b}{\phi} = \cls{(\cls{a}{\mu}+\cls{b}{\phi})}{\mu+\phi} =
    \cls{(\img{\iota_X}{\cls{a}{\mu}}\vee\img{\iota_Y}{\cls{b}{\phi}})}{\mu+\phi} \\
    = \cls{(\cls{\img{\iota_X}{a}}{\mu+\phi}\vee\cls{\img{\iota_Y}{b}}{\mu+\phi})}{\mu+\phi} =
    \cls{\cls{(\img{\iota_X}{a}\vee\img{\iota_Y}{b})}{\mu+\phi}}{\mu+\phi} =
    \cls{(a+b)}{\mu+\phi},
  \end{multline*}
  from idempotence of the closure operator \cite[see][Theorem
  3.1(a)]{2021-clos}, completing the proof.
\end{proof}

\msection{Sum of factorisations}
\label{sec:sum-of-fact}

The \fact{\mathsf{E}}{\mathsf{M}} factorisation of $f + g$ is already
achieved in Theorem \ref{sumisadmissible-equivalents}. When finite
sums of admissible subobjects is an admissible subobject,
factorisations also \emph{factor}.

\begin{Thm}
  \label{factsumcommute}
  In an extensive context with finite sum of admissible subobjects an
  admissible subobject, for any \Arr{f}{A}{X} and \Arr{g}{B}{Y}, with
  $f = \comp{m_X}{e_X}$ and $g = \comp{m_Y}{e_Y}$ the
  \fact{\mathsf{E}}{\mathsf{M}} of $f$ and $g$ respectively,
  $(f + g) = \comp{(m_X + m_Y)}{(e_X + e_Y)}$ is the
  \fact{\mathsf{E}}{\mathsf{M}} of $f + g$.
\end{Thm}

\begin{proof}
  \begin{figure}
    \centering
  \begin{equation*}
    % \resizebox{0.9\textwidth}{0.24cm}{
    \xymatrixcolsep{2.4em}
    \xymatrixrowsep{3em}
    \xymatrix{
      {A} \ar@{>->}[]!<1.2em,0em>;[rrr]^-{\iota_A} \ar[ddd]_f
      \ar@{-->}[drr]|{!\, e'_X}
      \ar@{->>}[ddr]^{e_X} & & & {A + B}
      \ar[ddd]|{f+g}
      \ar@{.>}[ddr]|(0.66){!\, e}
      \ar@{->>}[drr]|{e'} & & &
      {B} \ar[ddd]^g \ar@{->>}[ddr]^(0.66){e_Y}
      \ar@{>->}[]!<-1.2em,0em>;[lll]_-{\iota_B}
      \ar@{-->}[drr]|{!\, e'_Y} & & \\
      & & {T_X} \ar@{>->}[]!<-1.2em,-1.2em>;[dl]_-{m'_X}
      \ar[rrr]|(0.372){\hole}|(0.552){\hole}|(0.72){t_X} & & & {T}
      \ar@{>->}[]!<-1.2em,-1.2em>;[dl]|-{m'}
      \ar[rrr]|(0.324){\hole}|(0.48){\hole}|(0.72){t_Y} & & & {T_Y}
      \ar@{>->}[]!<-1.2em,-1.2em>;[dl]^-{m'_Y} \\
      & {I_X} \ar@{>->}[]!<-1.2em,-1.2em>;[dl]^-{m_X}
      \ar@{>->}[]!<1.2em,0em>;[rrr]|{i_X}|(0.636){\hole} & & & {I}
      \ar@{>->}[]!<-1.2em,-1.2em>;[dl]|-{m_X + m_Y} & & & {I_Y}
      \ar@{>->}[]!<-1.2em,-1.2em>;[dl]^-{m_Y}
      \ar@{>->}[]!<-1.2em,0em>;[lll]|(0.282){\hole}|{i_Y} & \\
      {X} \ar@{>->}[]!<1.2em,0em>;[rrr]_{\iota_X} & & & {X + Y}
      & & &
      {Y} \ar@{>->}[]!<-1.2em,0em>;[lll]^-{\iota_Y} & 
    }
    % }.
  \end{equation*}    
    \caption{Factorisation of sum}
    \label{fig:factsumcommute-diag}
  \end{figure}
  Consider the diagram in Figure \ref{fig:factsumcommute-diag}, where 
  the coproducts $X + Y$ and $A + B$ are presented in the bottom and
  top rows. Further, the morphisms $f$ and $g$ from \Bb{A} are given.
  From hypothesis, $m_X + m_Y \in \Sub{X + Y}{\mathsf{M}}$ and using
  extensivity the squares forming the coproduct $I = I_X + I_Y$ are
  both pullback squares.  There exists the unique morphism \Arr{e}{A +
    B}{I} such that $\comp{e}{\iota_A} = \comp{i_X}{e_X}$ and
  $\comp{e}{\iota_B} = \comp{i_Y}{e_Y}$.  Thus $e = e_X + e_Y$, and using
  extensivity the squares are pullback squares.  Let
  $e = \comp{m'}{e'}$ be the \fact{\mathsf{E}}{\mathsf{M}} for $e$.
  The admissible subobject $m'$ is pulled back along the coproduct
  injections $i_X$, $i_Y$ to obtain the pullback squares and hence the
  admissible subobjects $T_X$ and $T_Y$ of $I_X$ and $I_Y$
  respectively. In particular, from extensivity, $T = T_X + T_Y$,
  $t_X$, $t_Y$ are the coproduct injections and $m' = m'_X + m'_Y$.
  Hence there exist the unique morphisms $e'_X$ and $e'_Y$ such that
  $\comp{m'_x}{e'_X} = e_X$, $\comp{t_X}{e'_X} = \comp{e'}{\iota_A}$,
  $\comp{m'_Y}{e'_Y} = e_Y$ and
  $\comp{t_Y}{e'_Y} = \comp{e'}{\iota_B}$.  Since \down{e_X}{m'_X},
  \down{e_Y}{m'_Y}, $m'_X$ and $m'_Y$ are isomorphisms.  Hence $m'$ is
  also an isomorphism, proving $e = \comp{m'}{e'} \in \mathsf{E}$,
  completing the proof.
\end{proof}

\begin{Cor}
  \label{factsumcommute2}
  In an extensive context with finite sum of admissible subobjects an
  admissible subobject, for any \Arr{f}{A}{X}, \Arr{g}{B}{Y},
  $m_A \in \Sub{A}{\mathsf{A}}$ and $m_B \in \Sub{B}{\mathsf{M}}$:
  \begin{equation*}
    \begin{aligned}
      \comp{(f + g)}{(m_A + m_B)} & = \comp{f}{m_A} + \comp{g}{m_B}, \\
      \img{f+g}{(m_A + m_B)} & = \img{f}{m_A} + \img{g}{m_B}, \\
      \rest{(f + g)}{(M_A + M_B)} & = \rest{f}{M_A} + \rest{g}{M_B}.
    \end{aligned}
  \end{equation*}
\end{Cor}

\begin{proof}
    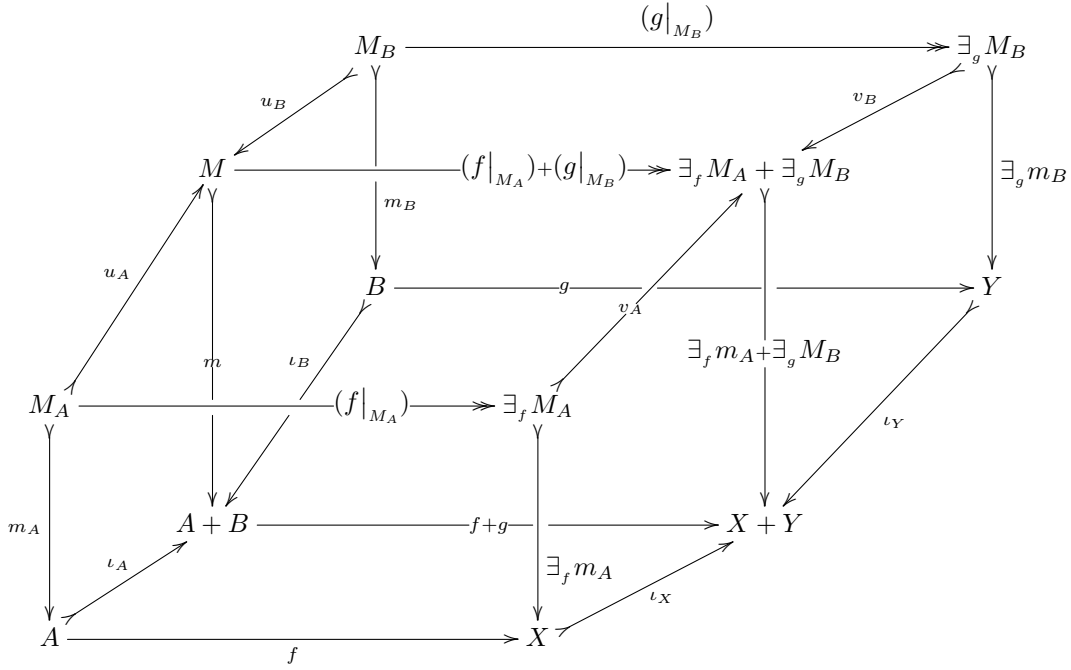
\begin{figure}
    \centering
  \begin{equation*}
    % \resizebox{36em}{!}{
    \xymatrixcolsep{3.36em}
    \xymatrixrowsep{3em}
    \xymatrix{
      & & {M_B} \ar@{>->}[]!<0pt,-1.2em>;[dd]|-{\hole}^(0.6){m_B}
      \ar@{->>}[rrr]^-{\rest{g}{M_B}}
      \ar@{>->}[]!<-1.2em,-1.2em>;[dl]_-{u_B}
      & & & {\img{g}{M_B}}
      \ar@{>->}[]!<0pt,-1.2em>;[dd]^-{\img{g}{m_B}}
      \ar@{>->}[]!<-1.56em,-0.96em>;[dl]_-{v_B} \\
      & {M} \ar@{>->}[]!<0em,-1.2em>;[ddd]|{m}|(0.636){\hole}
      \ar@{->>}[rrr]|(0.6){\rest{f}{M_A} +
        \rest{g}{M_B}}
      & & & {\img{f}{M_A} + \img{g}{M_B}}
      \ar@{>->}[]!<0pt,-1.2em>;[ddd]|-{\img{f}{m_A} + \img{g}{M_B}}
      & \\
      & & {B} \ar[rrr]|(0.306){g}|(0.444){\hole}|(0.638){\hole}
      \ar@{>->}[]!<-0.6em,-0.96em>;[ddl]_-(0.3){\iota_B}|(0.444){\hole}
      & & & {Y} \ar@{>->}[]!<-0.96em,-0.96em>;[ddl]^-{\iota_Y} \\
      {M_A} \ar@{>->}[]!<0pt,-1.2em>;[dd]_-{m_A}
      \ar@{>->}[]!<0.96em,0.96em>;[uur]^-{u_A}
      \ar@{->>}[rrr]|(0.66){\rest{f}{M_A}}
      & & & {\img{f}{M_A}} \ar@{>->}[]!<0.96em,1.02em>;[uur]|(0.336){v_A}
      \ar@{>->}[]!<0pt,-1.2em>;[dd]^(0.66){\img{f}{m_A}}& & \\
      & {A + B} 
      \ar[rrr]|{f + g}|(0.594){\hole}
      & & &
      {X + Y}  & \\
      {A} \ar@{>->}[]!<0.96em,0.96em>;[ur]^-{\iota_A}
      \ar[rrr]_-{f} & & &
      {X} \ar@{>->}[]!<1.2em,0.48em>;[ur]_-{\iota_X} & & \\    
    }
    % }
  \end{equation*}    
    \caption{Sum of factorisations}
    \label{fig:factsum-diag}
  \end{figure}

  Consider the diagram in Figure \ref{fig:factsum-diag}
  where $m_A \in \Sub{A}{\mathsf{M}}$ and
  $m_B \in \Sub{B}{\mathsf{M}}$, since
  \begin{equation*}
    \begin{aligned}
      \comp{f}{m_A} & = \comp{\img{f}{m_A}}{\rest{f}{M_A}}, \\
      \comp{g}{m_B} & = \comp{\img{g}{m_B}}{\rest{g}{M_B}}
    \end{aligned}
  \end{equation*}
  are the \fact{\mathsf{E}}{\mathsf{M}} for \comp{f}{m_A} and
  \comp{g}{m_B}, respectively. Using Theorem:
  \begin{equation*}\label{factsum-eq}\tag{$\dagger\star\dagger$}
    \begin{aligned}
      \comp{f}{m_A} + \comp{g}{m_B} & =
      \comp{\img{f}{m_A}}{\rest{f}{M_A}}
      + \comp{\img{g}{m_B}}{\rest{g}{M_B}} \\
      & = \comp{\bigl(\img{f}{m_A} +
        \img{g}{m_B}\bigr)}{\biggl(\rest{f}{M_A} +
        \rest{g}{M_B}\biggr)}
    \end{aligned}
  \end{equation*}
  is the \fact{\mathsf{E}}{\mathsf{M}} for
  $\comp{f}{m_A} + \comp{g}{m_B}$.  Hence, using extensivity, the
  vertical squares on the left and the right, the horizontal squares
  on the top and the bottom are all pullback squares. Hence the
  composites of the top and the vertical right hand squares are also
  pullbacks. Consequently, using extensivity again
  $\comp{(f + g)}{m} = \comp{f}{m_A} + \comp{g}{m_B}$.
\end{proof}

\msection{Pullback stability}
\label{sec:pb-closed-monic-E-morphisms}

The dense morphisms between finite sums is stable under pullbacks
along coproduct injections if and only if the closed embedding are
closed under finite sums (see Theorem
\ref{sumofclosedsubobjects}). Also, preneighbourhood morphisms from
\textsf{E} are dense \cite[see][Theorem 5.1(f)]{2021-clos}. Hence the
pullback along coproduct injections of any closed
\textsf{E}-monomorphism between finite sums is dense closed
monomorphism whenever finite sums of closed embeddings is a closed
embedding. The following theorem sharpens this observation.

\begin{Thm}
  \label{pbstabilityofclosedmonicformalepi}
  In an extensive context with finite sums of closed embeddings a
  closed embedding, the closed morphisms in
  $\mathsf{E} \cap \Mono{\Bb{A}}$ between finite sums are stable under
  pullbacks along coproduct injections.
\end{Thm}

\begin{proof}
  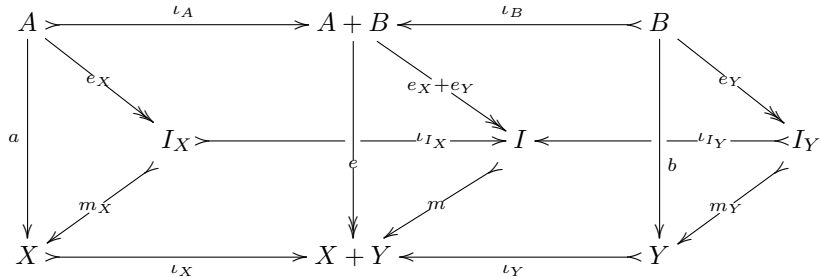
\begin{figure}[b]
    \centering
  \begin{equation*}
    % \resizebox{0.9\textwidth}{0.24cm}{
    \xymatrixcolsep{4em}
    \xymatrixrowsep{3em}
    \xymatrix{
      {A} \ar@{>->}[]!<1.2em,0em>;[rr]^-{\iota_A} \ar[dd]_-a
      \ar@{->>}[dr]|-{ e_X}
      & & {A + B} 
      \ar@{->>}[dd]|(0.6){ e}
      \ar@{->>}[dr]|-{ e_X + e_Y}
      & & {B} \ar[dd]^(0.6){b} \ar@{>->}[]!<-1.2em,0em>;[ll]_-{\iota_B}
      \ar@{->>}[dr]|-{ e_Y} \\
      & {I_X} \ar@{>->}[]!<-1.2em,-1.2em>;[dl]|-{ m_X}
      \ar@{>->}[]!<1.2em,0em>;[rr]|(0.462){\hole}|(0.72){\iota_{I_X}} & &
      {I} \ar@{>->}[]!<-1.2em,-1.2em>;[dl]|-{ m}
      & & {I_Y} \ar@{>->}[]!<-1.2em,-1.2em>;[dl]|-{ m_Y}
      \ar@{>->}[]!<-1.2em,0em>;[ll]|(0.24){\iota_{I_Y}}|(0.456){\hole} \\
      {X} \ar@{>->}[]!<1.2em,0em>;[rr]_-{\iota_X} & & {X + Y}
      & & {Y} \ar@{>->}[]!<-1.2em,0em>;[ll]^-{\iota_Y}
    }
    % }
  \end{equation*}
  \caption{Pullback stability of closed monomorphisms in \textsf{E}}
    \label{fig:pb-stab-closed-mono-in-E}
  \end{figure}

  Consider the diagram in Figure \ref{fig:pb-stab-closed-mono-in-E}{}
  where the top and bottom row are coproduct diagrams as shown and
  $e \in \mathsf{E} \cap \Mono{\Bb{A}}$ is a closed morphism.  Since
  $e \in \Mono{\Bb{A}}$ there exist unique $a \in \Sub{X}{}$ and
  $b \in \Sub{Y}{}$ such that $e = a + b$. Using extensivity the
  vertical squares are pullback squares. Since $e$ is closed and the
  coproduct injections continuous, the morphisms $a$ and $b$ are
  closed.  Let $a = \comp{m_X}{e_X}$ and $b = \comp{m_Y}{e_Y}$ be the
  \fact{\mathsf{E}}{\mathsf{M}} of $a$ and $b$ respectively.  Since
  $a$ is a closed morphism, for every $u \in \mc{\rest{\mu}{A}}$,
  $\cls{\img{a}{u}}{\mu} \leq \img{a}{u}$.  Since
  $m_X = \img{a}{\id{A}}$ and $\id{A} \in \mc{\rest{\mu}{A}}$,
  $\cls{m_X}{\mu} = \cls{\img{a}{\id{A}}}{\mu} \leq \img{a}{\id{A}} = m_X$
  implying $m_X$ to be closed. Similarly $m_Y$ is closed.  From
  hypothesis, $m = m_X + m_Y \in \mc{\mu+\phi}$ and the horizontal lower
  squares are pullback squares.  Since the top row is a coproduct,
  there exists the unique morphism $p = e_X + e_Y \in
  \mathsf{E}$. Therefore $\comp{m}{(e_X + e_Y)} = e$ is the
  \fact{\mathsf{E}}{\mathsf{M}} factorisation of $e \in
  \mathsf{E}$. Hence $m$ is an isomorphism.  From the horizontal
  pullback squares, $m_X$ and $m_Y$ are also isomorphisms. Hence
  $a, b \in \mathsf{E}$, completing the proof.
\end{proof}

%%% Local Variables:
%%% mode: latex
%%% TeX-master: "Extensivity"
%%% End:

%\input{CoprodInjectionsReflectZero}  %% section

\msection{Admissible subobjects of finite sums as biproduct}
\label{sec:sub-fin-sums-biprd}

Let \opair{X}{\mu} and \opair{Y}{\phi} be internal preneighbourhood spaces
in an extensive context \CAL{A}{.} Since each \Sub{X}{\mathsf{M}} a
distributive complete lattice \Fil{X} is distributive
\cite[see][Proposition 2.7(d)]{2020}, ensuring each \mc{\mu} is a
distributive lattice. Hence each \mc{\mu} can be considered as
$\vee$-semilattices, making them commutative monoids.

Assume \Arr{\iota_{X, Y}}{\Sub{X + Y}{\mathsf{M}}}{\Sub{X}{\mathsf{M}} \times
  \Sub{Y}{\mathsf{M}}} to be an isomorphism of partially ordered
sets. If $a, a' \in \Sub{X}{\mathsf{M}}$ and
$b, b' \in \Sub{Y}{\mathsf{M}}$ then using Theorem
\ref{sumisadmissible-equivalents}:
\begin{equation*}
  \begin{aligned}
    (a + b) \vee (a' + b') & = (\img{\iota_X}{a} \vee \img{\iota_Y}{b})
    \vee (\img{\iota_X}{a'} \vee \img{\iota_Y}{b'}) \\
    & = (\img{\iota_X}{a} \vee \img{\iota_X}{a'}) \vee (\img{\iota_Y}{b}
    \vee \img{\iota_Y}{b'}) \\
    & = \img{\iota_X}{(a \vee a')} \vee \img{\iota_Y}{(b \vee b')} \\
    & = (a \vee a') + (b \vee b')
  \end{aligned}
\end{equation*}
implies
$\finv{\iota_X}{\bigl((a + b) \vee (a' + b')\bigr)} = a \vee a'$ and
$\finv{\iota_Y}{\bigl((a + b) \vee (a' + b')\bigr)} = b \vee b'$. As a
consequence, the diagram
\begin{equation}
  \label{subofsum=biproductofsummands-eq1}
  \xymatrixcolsep{1.2cm}
  \xymatrix{
    {\Sub{X}{\mathsf{M}}} \ar@<1ex>[r]^-{\img{\iota_X}{}}
    \ar@<-1ex>@{<-}[r]_-{\finv{\iota_X}{}} \ar@{}[r]|-{\bot} &
    {\Sub{X + Y}{\mathsf{M}}}
    \ar@{<-}@<1ex>[r]^-{\img{\iota_Y}{}}
    \ar@<-1ex>[r]_-{\finv{\iota_Y}{}} \ar@{}[r]|-{\bot} &
    {\Sub{Y}{\mathsf{M}}}               
  }
\end{equation}
is a biproduct in the category \SemLat{\vee} of $\vee$-semilattices and
$\vee$-semilattice homomorphisms.  The case for closed subobjects is
similar --- if \Arr{\iota_{X, Y}}{\mc{\mu + \phi}} {\mc{\mu} \times \mc{\phi}} is an
isomorphism of partially ordered sets, then
$\cls{(\img{\iota_X}{-} \vee \img{\iota_Y}{-})}{\mu + \phi} = \inv{\iota_{X,Y}} = +_{X, Y}
= \img{\iota_X}{} \vee \img{\iota_Y}{}$.

\begin{Thm}
  \label{subobjectofsumisbiproduct}
  \index[res]{Theorem: subobjects of sum as biproduct} Let
  \opair{X}{\mu} and \opair{Y}{\phi} be internal preneighbourhood spaces in
  an extensive context \CAL{A}.  Let $\CAL{K}(-)$ stand for either
  \Sub{-}{\mathsf{M}} or \mc{-}.

  If \Arr{\iota_{X, Y}}{\CAL{K}(X + Y)}{\CAL{K}(X) \times \CAL{K}(Y)} be an
  isomorphism of partially ordered sets then
  \begin{equation*}
    \xymatrixcolsep{1.2cm}
    \xymatrix{
      {\CAL{K}(X)} \ar@<1ex>[r]^-{\img{\iota_X}{}}
      \ar@<-1ex>@{<-}[r]_-{\finv{\iota_X}{}} \ar@{}[r]|-{\bot} &
      {\CAL{K}(X + Y)}
      \ar@{<-}@<1ex>[r]^-{\img{\iota_Y}{}}
      \ar@<-1ex>[r]_-{\finv{\iota_Y}{}} \ar@{}[r]|-{\bot} &
      {\CAL{K}(Y)}               
    }
  \end{equation*}
  is a biproduct in \SemLat{\vee}.
\end{Thm}

\begin{rem}
  Thus, in context of Theorem, since $\CAL{K}(X + Y)$ is a biproduct,
  the $\vee$-semilattice homomorphisms \Arr{f}{\CAL{K}(X + Y)}{\CAL{K}(P
    + Q)} correspond to matrices
  \begin{equation*}
    \begin{pmatrix}
      \Arr{f_{X,P}}{\CAL{K}(X)}{\CAL{K}(P)} &
      \Arr{f_{X,Q}}{\CAL{K}(X)}{\CAL{K}(Q)} \\
      \Arr{f_{Y,P}}{\CAL{K}(Y)}{\CAL{K}(P)} &
      \Arr{f_{Y,Q}}{\CAL{K}(Y)}{\CAL{K}(Q)}
    \end{pmatrix}
  \end{equation*}
  natural in $X$, $Y$, $P$ and $Q$.
\end{rem}

\begin{rem}
  Theorem further states: the $\vee$-semilattice $\CAL{K}(X + Y)$ is a
  biproduct of the $\vee$-semilattices $\CAL{K}(X)$ and $\CAL{K}(Y)$ if
  and only if $\CAL{K}(X + Y)$ is isomorphic to
  $\CAL{K}(X) \times \CAL{K}(Y)$ as partially ordered sets.
\end{rem}

\msection{Sum of proper morphisms}
\label{sec:sum-proper}

Recall: an internal preneighbourhood space \opair{X}{\mu}{} is compact
if the unique preneighbourhood morphism
\Arr{\terma{X}}{\opair{X}{\mu}}{\opair{\termo}{\nabla_{\termo}}}{} is proper
\cite[see][Remark (K)]{2021-clos}, and \Comp{\Bb{A}}{} is the full
subcategory of \pNHD{\Bb{A}}{} consisting of all compact spaces
\cite[see][\S6.2{}]{2021-clos}.

\begin{Thm}
  \label{thm:sum-proper}{}
  In an extensive context \CAL{A} in which finite sum of closed
  morphisms is closed, a sum of proper morphisms is proper and
  \Comp{\Bb{A}}{} is closed under finite sums if and only if for each
  internal preneighbourhood space \opair{X}{\mu}, both the morphisms
  \begin{equation*}
    \xymatrixcolsep{6em}
    \xymatrix{
      {\opair{\inito}{\uparrow_{\inito}}} \ar[r]^-{\inita{X}} & {\opair{X}{\mu}}
      \ar@{<-}[r]^-{\Opair{\id{X}}{\id{X}}} & {\opair{X+X}{\mu+\mu}}
    } 
  \end{equation*}
  are proper.
\end{Thm}

\begin{proof}
  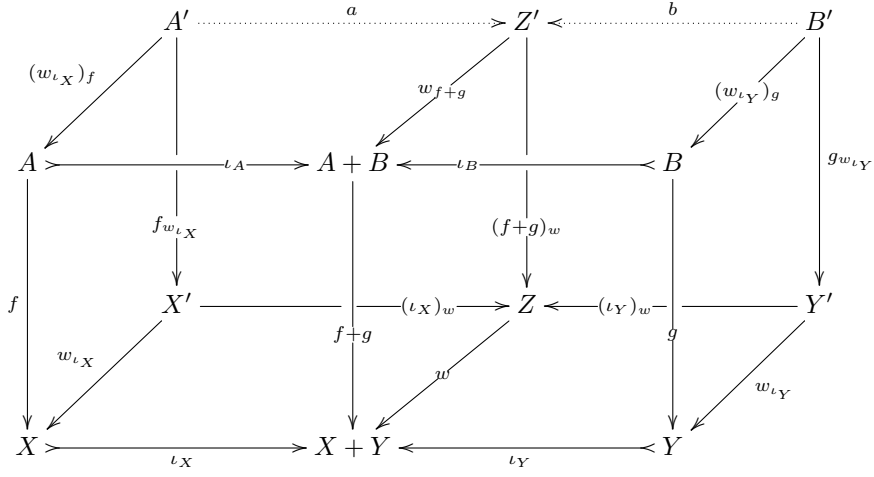
\begin{figure}
    \centering
    \begin{equation*}
      \xymatrixcolsep{4em}
      \xymatrixrowsep{4em}
      \xymatrix{
        & {A'} \ar[dd]|(0.492){\hole}|(0.72){f_{w_{\iota_X}}} \ar@{.>}[rr]^-a
        \ar[dl]_-{(w_{\iota_X})_{f}}
        & & {Z'} \ar[dd]|(0.492){\hole}|(0.72){(f+g)_{w}} \ar@{<.}[rr]^-b
        \ar[dl]|-{w_{f+g}}
        & & {B'} \ar[dd]^-{g_{w_{\iota_Y}}} \ar[dl]|-{(w_{\iota_Y})_g} \\
        {A} \ar@{>->}[]!<1.2em,0em>;[rr]|(0.6){\iota_A} \ar[dd]_-f
        & & {A + B}  \ar[dd]|(0.6){f+g}
        & & {B} \ar[dd]|(0.6){g}
        \ar@{>->}[]!<-1.2em,0em>;[ll]|(0.6){\iota_B} \\
        & {X'} \ar[rr]|(0.492){\hole}|(0.72){(\iota_X)_{w}} \ar[dl]_-{w_{\iota_X}}
        & & {Z} \ar@{<-}[rr]|(0.336){(\iota_Y)_{w}}|(0.504){\hole}
        \ar[dl]|-{w}
        & & {Y'} \ar[dl]^-{w_{\iota_Y}} \\
        {X} \ar@{>->}[]!<1.2em,0em>;[rr]_-{\iota_X}
        & & {X + Y} 
        & & {Y} \ar@{>->}[]!<-1.2em,0em>;[ll]^-{\iota_Y}
      }
    \end{equation*}    
    \caption{Sum of proper morphisms}
    \label{fig:sum-proper}
  \end{figure}

  Given the proper morphisms \Arr{f}{\opair{A}{\alpha}}{\opair{X}{\mu}}{,}
  \Arr{g}{\opair{B}{\beta}}{\opair{Y}{\phi}} and preneighbourhood
  morphism \Arr{w}{\opair{Z}{\psi}}{\opair{X+Y}{\mu+\phi}}{} consider the
  diagram in Figure \ref{fig:sum-proper} where the middle vertical is
  the pullback of $f + g$ along $w$. Extensivity implies the front two
  vertical squares are pullback squares. The morphism $w$ is pulled
  back along the coproduct injections $\iota_X$ and $\iota_Y$, yielding from
  extensivity $Z = X' + Y'$ and $w = w_{\iota_X} + w_{\iota_Y}$. The morphisms
  $f$ and $g$ are pulled back along $w_{\iota_X}$ and $w_{\iota_Y}$, yielding
  from extensivity again the left hand and right hand vertical
  squares. Since
  $\comp{w}{\comp{(\iota_X)_w}{f_{w_{\iota_X}}}} =
  \comp{(f+g)}{\comp{\iota_A}{(w_{\iota_X})_f}}$ and
  $\comp{w}{\comp{(\iota_Y)_w}{g_{w_{\iota_Y}}}} =
  \comp{(f+g)}{\comp{\iota_B}{(w_{\iota_Y})_g}}$, from the middle vertical
  pullback square there exist unique morphisms $a$ and $b$ which makes
  the whole diagram to commute. Hence from properties of pullback
  squares all the squares in the diagram are pullback
  squares. Consequently, from extensivity, $Z' = A' + B'$,
  $(f + g)_{w} = f_{w_{\iota_X}} + g_{w_{\iota_Y}}$. Since $f$ and $g$ are
  proper morphisms, $f_{w_{\iota_X}}$ and $g_{w_{\iota_Y}}$ are closed
  morphisms. Since sum of closed morphisms is closed, $(f+g)_{w}$ is a
  closed morphism, proving $f + g$ is proper. This proves the first
  part of the statement. Using properties of proper morphisms
  \cite[see][Theorem 6.1(a) \& Remark (K)]{2021-clos}
  $\opair{\inito}{\uparrow_{\inito}}$ is compact if and only if for any
  internal preneighbourhood space \opair{X}{\mu}{} the projection
  \Arr{p_2}{\opair{\inito\times{X}}{\uparrow_{\inito}\times\mu}}{\opair{X}{\mu}} is proper,
  i.e., \Arr{\inita{X}}{\opair{\inito}{\uparrow_{\inito}}}{\opair{X}{\mu}}{} is
  proper, since $\inito\times{X} \approx \inito$. Similarly,
  $\opair{\termo+\termo}{\nabla_{\termo}+\nabla_{\termo}}$ is compact if and
  only if for any internal preneighbourhood space \opair{X}{\mu}{,} the
  second projection \Arr{p_2}{\opair{(\termo+\termo)\times{X}}
    {(\nabla_{\termo}+\nabla_{\termo})\times\mu}}{\opair{X}{\mu}}{} is proper. Since
  every lextensive category is distributive \cite[see][Proposition
  4.5]{CarboniLackWalters1993}, this is equivalent to the morphism
  \Arr{\Opair{\id{X}}{\id{X}}}{\opair{X+X}{\mu+\mu}}{\opair{X}{\mu}} being
  proper. This proves the \emph{only if} side of the second part; the
  \emph{if} side follows from the first part trivially, completing the
  proof.
\end{proof}

Thus, continuing from \cite[][Theorem 6.2(c)]{2021-clos}, in an extensive context $\CAL{A} = (\Bb{A}, \mathsf{E}, \mathsf{M})$ with admissible subobjects closed under finite sums, coproduct injections closed embeddings and each \inita{X}{} and codiagonal \Opair{\id{X}}{\id{X}}{} closed, the cateory \Comp{\Bb{A}}{} of compact internal preneighbourhood spaces is finitely productive, closed hereditary (i.e., every closed embedding to a compact space is compact) and closed under finite sums.

\msection{Sum of separated morphisms}
\label{sec:sum-of-separated}

Recall: an internal preneighbourhood space \opair{X}{\mu}{} is Hausdorff
if the unique preneighbourhood morphism
\Arr{\terma{X}}{\opair{X}{\mu}}{\opair{\termo}{\nabla_{\termo}}}{} is
separated and \Int{\Haus}{\Bb{A}}{} is the full subcategory of
internal Hausdorff spaces.

\begin{Thm}
  \label{thm:sum-of-separated}
  In an extensive context with finite sum of closed morphisms closed,
  a sum of separated morphisms is separated and \Int{\Haus}{\Bb{A}}{}
  is closed under finite sums if and only if for
  \opair{\termo+\termo}{\nabla_{\termo}+\nabla_{\termo}}{} is Hausdorff.
\end{Thm}

\begin{proof}
  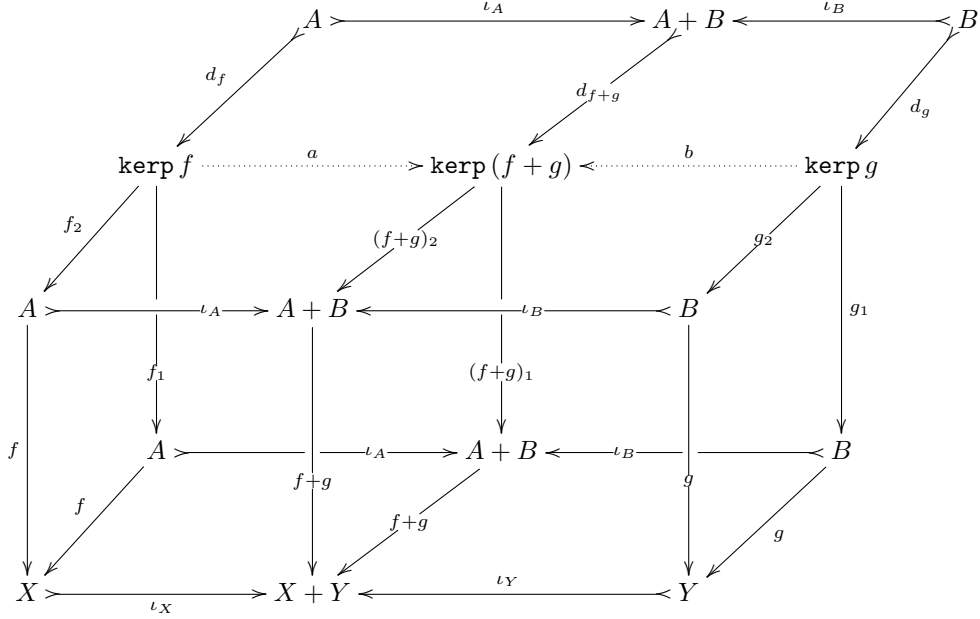
\begin{figure}
    \centering
    \begin{equation*}
      \xymatrixcolsep{2.4em}
      \xymatrixrowsep{4em}
      \xymatrix{
        & & {A} \ar@{>->}[]!<1.2em,0em>;[rr]^-{\iota_A}
        \ar@{>->}[]!<-0.78em,-0.72em>;[dl]_-{d_f}
        & & {A + B} \ar@{>->}[]!<-0.78em,-0.72em>;[dl]|-{d_{f+g}}
        & & {B} \ar@{>->}[]!<-0.78em,-0.72em>;[dl]^-{d_g}
        \ar@{>->}[]!<-1.2em,0em>;[ll]_-{\iota_B} \\
        & {\Kerp{f}}
        \ar[dd]|(0.492){\hole}|(0.72){f_1} \ar@{.>}[rr]^-a
        \ar[dl]_-{f_{2}}
        & & {\Kerp{(f+g)}}
        \ar[dd]|(0.492){\hole}|(0.72){(f+g)_{1}} \ar@{<.}[rr]^-b
        \ar[dl]|-{(f+g)_{2}}
        & & {\Kerp{g}} \ar[dd]^-{g_{1}} \ar[dl]|-{g_{2}} \\
        {A} \ar@{>->}[]!<1.2em,0em>;[rr]|(0.6){\iota_A} \ar[dd]_-f
        & & {A + B} \ar[dd]|(0.6){f+g}
        & & {B} \ar[dd]|(0.6){g}
        \ar@{>->}[]!<-1.2em,0em>;[ll]|(0.36){\iota_B} \\
        & {A} \ar@{>->}[]!<1.2em,0em>;[rr]|(0.396){\hole}|(0.6){\iota_{A}}
        \ar[dl]_-{f}
        & & {A+B}
        \ar[dl]|-{f+g}
        & & {B} \ar[dl]^-{g}
        \ar@{>->}[]!<-1.2em,0em>;[ll]|(0.388){\hole}|(0.6){\iota_B} \\
        {X} \ar@{>->}[]!<1.2em,0em>;[rr]_-{\iota_X}
        & & {X + Y} 
        & & {Y} \ar@{>->}[]!<-1.2em,0em>;[ll]_-{\iota_Y}
      }
    \end{equation*}    
    \caption{Sum of separated morphisms}
    \label{fig:sum-of-separated}
  \end{figure}

  Given the separated morphisms \Arr{f}{\opair{A}{\alpha}}{\opair{X}{\mu}},
  \Arr{g}{\opair{B}{\beta}}{\opair{Y}{\phi}}{,} consider the diagram in
  Figure \ref{fig:sum-of-separated} where the three vertical squares
  are the kernel pairs of $f$, $f+g$ and $g$, lower horizontal and
  front vertical squares are the sums, which by extensivity are also
  pullback squares. Since
  $\comp{(f+g)}{\comp{\iota_A}{f_1}} = \comp{(f+g)}{\comp{\iota_A}{f_2}}$ and
  $\comp{(f+g)}{\comp{\iota_B}{g_1}} = \comp{(f+g)}{\comp{\iota_B}{g_2}}$
  there exist unique morphisms $a$ and $b$, making the whole diagram
  to commute. An use of properties of pullbacks ensure every square in
  the front cube are all pullbacks. Hence, using extensivity,
  $(f+g)_{i} = f_{i} + g_{i}$ ($i = 1, 2$),
  $\Kerp{(f+g)} = \Kerp{f} + \Kerp{g}$ and the morphisms $a$, $b$ are
  coproduct injections.  Clearly, $d_{f}$, $d_{f+g}$ and $d_{g}$ are
  equalisers of the kernel pairs \opair{f_1}{f_2}{,}
  \opair{(f+g)_1}{(f+g)_2}{} and \opair{g_1}{g_2}{.} Using extensivity
  again, $d_{f}$ (or, $d_{g}$) is the pullback of $d_{f+g}$ along $a$
  (or, $b$), implying $d_{f+h} = d_{f} + d_{g}$. Since $f$ and $g$ are
  separated, $d_{f}$, $d_{g}$ are both proper, so that an use of
  Theorem \ref{thm:sum-proper}{} implies $d_{f+g}$ is proper, i.e.,
  $f + g$ is separated. The \emph{only if} side of the second part if
  trivial; conversely, if
  \opair{\termo+\termo}{\nabla_{\termo}+\nabla_{\termo}}{} is Hausdorff, then
  from the first part for any two Hausdorff spaces
  $\terma{X}+\terma{Y}$ is separated,
  $\terma{X+Y} = \comp{\terma{\termo+\termo}}{(\terma{X}+\terma{Y})}$
  being composite of two separated morphisms is separated
  \cite[see][Theorem 7.1]{2021-clos} i.e., $X + Y$ is Hausdorff,
  completing the proof.
\end{proof}

Thus, continuing from \cite[][Corollary 7.2]{2021-clos}, in an extensive context $\CAL{A} = (\Bb{A}, \mathsf{E}, \mathsf{M})$ with admissible subobjects closed under finite sums, coproduct injections closed embeddings and \opair{\termo+\termo}{\nabla_{\termo}+\nabla_{\termo}}{} an internal Hausdorff preneighbourhood space, the cateory \Int{\Haus}{\Bb{A}} of Hausdorff internal preneighbourhood spaces is finitely complete, hereditary (i.e., every subobject of an internal Hausdorff preneighbourhood space is Hausdorff), closed under images of preneighbourhood morphisms stably in \textsf{E} and is closed under finite sums.

%backmatter

\section*{Acknowledgements}\label{acknowledgements}
\,\newline
I am gratefully indebted:
\begin{enumerate}
\item To the funding received from the {\em European Union Horizon
    2020 MCSA Irses project 731143} and is thankfully acknowledged.

\item To the useful suggestions received from Amartya Goswami, Themba
  Dube and George Janelidze who painstakingly went through
  the first draft of this paper. The errors that remain are entirely
  my misconceptions which needs to be purged with time and
  forbearance.
\end{enumerate}

\printbibliography

% \printindex[def] \printindex[res] \printindex[sym]

%%% insert the right number of ../'s from the current directory
% \bibliographystyle{ams}
% \bibliography{../../../tex/essentials/bib-algebra,../../../tex/essentials/bib-category,../../../tex/essentials/bib-frames,../../../tex/essentials/bib-topology}

\end{document}

%%% Local Variables:
%%% mode: latex
%%% TeX-master: t
%%% End: